\documentclass[leqno,12pt]{article}

\usepackage{amssymb}
\usepackage{euscript}
\usepackage[dvips]{graphicx}
\usepackage{flafter}
\usepackage{pstricks}

\usepackage{verbatim}

\usepackage{amsmath}

\usepackage{color}

\usepackage{mathrsfs} 

\usepackage{amsthm}

\evensidemargin\oddsidemargin

\begin{document}

\baselineskip=18pt
\setcounter{page}{1}

\renewcommand{\theequation}{\thesection.\arabic{equation}}
\newtheorem{theorem}{Theorem}[section]
\newtheorem{lemma}[theorem]{Lemma}
\newtheorem{definition}[theorem]{Definition}
\newtheorem{proposition}[theorem]{Proposition}
\newtheorem{corollary}[theorem]{Corollary}
\newtheorem{fact}[theorem]{Fact}
\newtheorem{problem}[theorem]{Problem}
\newtheorem{conjecture}[theorem]{Conjecture}
\newtheorem{claim}[theorem]{Claim}

\theoremstyle{definition} 
\newtheorem{remark}[theorem]{Remark}

\newcommand{\eqnsection}{
\renewcommand{\theequation}{\thesection.\arabic{equation}}
    \makeatletter
    \csname  @addtoreset\endcsname{equation}{section}
    \makeatother}
\eqnsection


\def\r{{\mathbb R}}
\def\e{{\mathbb E}}
\def\p{{\mathbb P}}
\def\P{{\bf P}}
\def\E{{\bf E}}
\def\Q{{\bf Q}}
\def\z{{\mathbb Z}}
\def\N{{\mathbb N}}
\def\T{{\mathbb T}}
\def\G{G}

\def\ee{\mathrm{e}}
\def\d{\, \mathrm{d}}



\vglue50pt

\centerline{\Large\bf The most visited sites of biased random walks on trees}

{
\let\thefootnote\relax\footnotetext{\scriptsize Partly supported by ANR project MEMEMO2 (2010-BLAN-0125).}
}

\bigskip
\bigskip

\centerline{by}

\medskip

\centerline{Yueyun Hu\footnote{\scriptsize LAGA, Universit\'e Paris XIII, 99 avenue J-B Cl\'ement, F-93430 Villetaneuse, France, {\tt yueyun@math.univ-paris13.fr}} and Zhan Shi\footnote{\scriptsize LPMA, Universit\'e Paris VI, 4 place Jussieu, F-75252 Paris Cedex 05, France, {\tt zhan.shi@upmc.fr}}}

\medskip

\centerline{\it Universit\'e Paris XIII \& Universit\'e Paris VI}



\bigskip
\bigskip
\bigskip

{\leftskip=2truecm \rightskip=2truecm \baselineskip=15pt \small

\noindent{\slshape\bfseries Summary.} We consider the slow movement of randomly biased random walk $(X_n)$ on a supercritical Galton--Watson tree, and are interested in the sites on the tree that are most visited by the biased random walk. Our main result implies tightness of the distributions of the most visited sites under the annealed measure. This is in contrast with the one-dimensional case, and provides, to the best of our knowledge, the first non-trivial example of null recurrent random walk whose most visited sites are not transient, a question originally raised by Erd\H os and R\'ev\'esz~\cite{erdos-revesz} for simple symmetric random walk on the line. 

\bigskip

\noindent{\slshape\bfseries Keywords.} Biased random walk on the Galton--Watson tree, branching random walk, local time, most visited site.

\bigskip

\noindent{\slshape\bfseries 2010 Mathematics Subject
Classification.} 60J80, 60G50, 60K37.

} 

\bigskip
\bigskip

\section{Introduction}
   \label{s:intro}

$\phantom{aob}$We consider a (randomly) biased random walk $(X_n)$ on a  supercritical Galton--Watson tree    $\T$,  rooted at $\varnothing$.  The random biases are represented by $\omega := (\omega(x), \, x\in \T \backslash \{ \varnothing\} )$, a family of random vectors; for each vertex $x\in \T$, $\omega(x) := (\omega(x, \, y), \, y\in \T)$ is such that $\omega(x, \, y) \ge 0$ for all $y\in \T$ and that $\sum_{y\in \T} \omega(x, \, y) =1$.  For any vertex $x\in \T \backslash \{ \varnothing\}$, let ${\buildrel \leftarrow \over x}$ be its parent. For the sake of presentation, we modify the values of $\omega(\varnothing, \, x)$ for $x$ with ${\buildrel \leftarrow \over x} = \varnothing$, and add a special vertex, denoted by ${\buildrel \leftarrow \over \varnothing}$, which is considered as the parent of $\varnothing$, such that $\omega(\varnothing, \, {\buildrel \leftarrow \over \varnothing}) + \sum_{x: \, {\buildrel \leftarrow \over x} = \varnothing} \omega(\varnothing, \, x) =1$. The vertex ${\buildrel \leftarrow \over \varnothing}$ is, however, {\it not} regarded as a vertex of $\T$; so, for example, $\sum_{x\in \T} f(x)$ does not contain the term $f({\buildrel \leftarrow \over \varnothing})$. 

Assume that for each pair of vertices $x$ and $y$ in $\T\cup\{{\buildrel \leftarrow \over \varnothing}\}$, $\omega(x, \, y)>0$ if and only if $y \sim x$, where by $x\sim y$ we mean that  $x$ is either a child, or the parent, of $y$.   Moreover, we define  $\omega({\buildrel \leftarrow \over \varnothing}, \varnothing):=1$.

Given $\omega$,   the biased walk $(X_n, \, n\ge 0)$ is  a Markov chain taking values on $\T\cup\{{\buildrel \leftarrow \over \varnothing}\}$, started at $X_0 = \varnothing$, whose transition probabilities are
$$
P_\omega \{ X_{n+1} = y \, | \, X_n =x \} = \omega(x, \, y).
$$

\noindent The probability $P_\omega$ is often referred to as the quenched probability. We also consider the annealed probability $\p(\, \cdot \, ) := \int P_\omega( \, \cdot \, ) \, \P(\! \d \omega)$, where $\P$ denotes the probability with respect to the environment $(\omega, \T)$.

There is an active literature on randomly biased walks on Galton-Watson trees; see, for example, a large list of references in \cite{yzenergy}. In this paper, we restrict our attention to a regime of {\it slow movement} of the walk in the recurrent case.
 
Clearly, the movement of the biased random walk $(X_n)$ is determined by the law of the random environment $\omega$. We assume that $(\omega(x, \, y), \, y\sim x)$ for $x\in \T $, are i.i.d.\ random vectors. It is convenient to view $(\omega, \T)$ as a marked tree (in the sense of Neveu \cite{neveu}).
  
The influence of the random environment is quantified by means of the random {\bf potential} process $(V(x), \, x\in \T)$, defined by $V(\varnothing):=0$ and
\begin{equation}
    V(x) 
    := 
    -
    \sum_{y\in \, ]\!] \varnothing,\, x]\!]}
    \log \,
    \frac{\omega({\buildrel \leftarrow \over y},
    \, y)}{\omega({\buildrel \leftarrow \over y}, \,
    {\buildrel \Leftarrow \over y})},
    \qquad x\in \T\backslash\{ \varnothing\} \, ,
    \label{V}
\end{equation}

\noindent where ${\buildrel \Leftarrow \over y}$ is the parent of ${\buildrel \leftarrow \over y}$, and $\, ]\!] \varnothing, \, x]\!] := [\![ \varnothing, \, x]\!] \backslash \{ \varnothing\}$, with $[\![ \varnothing, \, x]\!]$ denoting the set of vertices (including $x$ and $\varnothing$) on the unique shortest path connecting $\varnothing$ to $x$. There exists an obvious bijection between the random environment $\omega$ and the random potential $V$.

For any $x\in \T$, let $|x|$ denote its generation. Throughout the paper, we assume   
\begin{equation}
    \E \Big( \sum_{x: \, |x|=1} \ee^{-V(x)} \Big) =1,
    \qquad
    \E \Big( \sum_{x: \, |x|=1} V(x)\, \ee^{-V(x)} \Big) =0 \, .
    \label{cond-hab}
\end{equation}

\noindent We also assume that the following integrability condition is fulfilled: there exists $\delta>0$ such that 
\begin{equation}
    \E \Big( \sum_{x: \, |x|=1} \ee^{-(1+\delta)V(x)} \Big)
    +
    \E \Big( \sum_{x: \, |x|=1} \ee^{ \delta V(x)} \Big) 
    +
    \E \Big[ \Big( \sum_{x: \, |x|=1} 1\Big)^{1+\delta} \, \Big] <\infty \, .
    \label{integrability-assumption}
\end{equation}

The random potential $(V(x), \, x\in \T)$ is a branching random walk as in Biggins~\cite{biggins77}; as such, (\ref{cond-hab}) corresponds to the ``boundary case" (Biggins and Kyprianou~\cite{biggins-kyprianou05}). It is known that, under some additional integrability assumptions that are weaker than \eqref{integrability-assumption}, the branching random walk in the boundary case possesses some deep universality properties, see \cite{stf} for references.

Under $(\ref{cond-hab})$ and $(\ref{integrability-assumption})$, the biased walk $(X_n)$ is null recurrent (Lyons and Pemantle~\cite{lyons-pemantle},  Menshikov and Petritis~\cite{menshikov-petritis}, Faraud~\cite{faraud}), such that upon the system's survival,
\begin{eqnarray}
    \frac{|X_n|}{(\log n)^2} 
 &\; \buildrel \mathrm{law} \over \longrightarrow \; &
    X_\infty\, , 
    \label{claw}
    \\
	\frac{1}{(\log n)^3} \, \max_{0\le i \le n}|X_i| 
 &\; \to \;&
    c_1 \; \qquad \hbox{\rm a.s.},
	\label{ascv}
\end{eqnarray} 

\noindent where $X_\infty$ is non-degenerate taking values in $(0, \, \infty)$, and $c_1$ denotes a positive constant: both $X_\infty$ and $c_1$ are explicitly known, see \cite{yzlocaltree} and \cite{gyzbiased}, respectively.

For any vertex $x\in \T$, let us define
$$
L_n(x)
:=
\sum_{i=1}^n \, {\bf 1}_{\{ X_i =x\} } \, ,
\qquad
n\ge 1 \, ,
$$

\noindent which is the (site) local time of the biased walk at $x$. Consider, for any $n\ge 1$, the non-empty random set
\begin{equation}
    \mathscr{A}_n
    :=
    \Big\{ x\in \T: \, L_n(x) = \max_{y\in \T}L_n(y)\Big\} .
    \label{ensemble_points_favoris}
\end{equation}

\noindent In words, $\mathscr{A}_n$ is the set of the most visited sites (or: favourite sites) at time $n$. The study of favourite sites was initiated by Erd\H os and R\'ev\'esz~\cite{erdos-revesz} for the symmetric Bernoulli random walk on the line (see a list of ten open problems presented in Chapter 11 of the book of R\'ev\'esz~\cite{revesz}). In particular, for the symmetric Bernoulli random walk on $\z$, Erd\H os and R\'ev\'esz~\cite{erdos-revesz} conjectured: (a) tightness for the family of most visited sites, and (b) the cardinality of the set of most visited sites being eventually bounded by $2$. Conjecture (b) was partially proved by T\'oth~\cite{toth}, and is believed to be true by many. On the other hand, Conjecture (a) was disproved by Bass and Griffin~\cite{bass-griffin}: as a matter of fact, $\inf \{ |x|, \; x\in \mathscr{A}_n\} \to \infty$ almost surely for the one-dimensional Bernoulli walk. Later, we proved in \cite{fav-rwre} that it was also the case for Sinai's one-dimensional random walk in random environment. The present paper is devoted to studying both questions for biased walks on trees; our answer is as follows. 

\bigskip

\noindent {\bf Corollary \ref{c:main}.} 
{\it Assume $(\ref{cond-hab})$ and $(\ref{integrability-assumption})$.  There exists a finite non-empty set $\mathscr{U}_{\min}$, defined in $(\ref{Umin})$ and depending only on the environment, such that 
 $$
 \lim_{n\to \infty} \, \p( \mathscr{A}_n \subset \mathscr{U}_{\min} \, | \, \hbox{non-extinction}) 
 =
 1 \, .
 $$

 \noindent In particular, the family of most visited sites is tight under $\p$.
} 

\bigskip

So, concerning the tightness question for most visited sites, biased walks on trees behave very differently from recurrent one-dimensional nearest-neighbour random walks (whether the environment is random or deterministic). To the best of our knowledge, this is the first non-trivial example of null recurrent Markov chain whose most visited sites are tight.

In the next section, we give a precise statement of the main result of this paper, Theorem \ref{t:main}.


 
\section{Statement of results}


$\phantom{aob}$Let us define a symmetrized version of the potential: 
\begin{equation}
    U(x) 
    := 
    V(x) - \log (\frac{1}{\omega(x, \, {\buildrel \leftarrow \over x})}) \, ,
    \qquad
    x\in \T \, .
    \label{U}
\end{equation}

\noindent Note that 
\begin{equation}
    \ee^{-U(x)}
    = 
    \frac{1}{\omega(x, \, {\buildrel \leftarrow \over x})} \, \ee^{-V(x)}  
    =
    \ee^{-V(x)} + \sum_{y\in \T: \, {\buildrel \leftarrow \over y} =x} \ee^{-V(y)} ,
    \qquad
    x\in \T \, .
    \label{Theta}
\end{equation}


It is known (Biggins~\cite{biggins-mart-cvg}, Lyons~\cite{lyons}) that under assumption (\ref{cond-hab}), 
\begin{equation}
    \inf_{x: \, |x|=n} U(x) 
    \to 
    \infty ,
    \qquad
    \hbox{\rm $\P^*$-a.s.},
    \label{U->infty}
\end{equation}

\noindent where here and in the sequel, 
\begin{eqnarray*}
    \P^*(\, \cdot \, )
 &:=&\P( \, \cdot \, | \, \hbox{non-extinction}) ,
    \\
    \p^*(\, \cdot \, )
 &:=&\p( \, \cdot \, | \, \hbox{non-extinction}) \, .
\end{eqnarray*}

Define the {\it derivative martingale}
\begin{equation}
    D_n 
  :=  \sum_{x: \, |x|=n} V(x) \ee^{-V(x)}, 
    \qquad
    n\ge 0 \, .
    \label{Dn}
\end{equation}

\noindent It is known (Biggins and Kyprianou~\cite{biggins-kyprianou04}, A\"id\'ekon~\cite{elie-min}, Chen~\cite{chen}) that (\ref{integrability-assumption}) implies that $D_n$ converges $\P$-a.s.\ to a limit, denoted by $D_\infty$, and that 
$$
D_\infty >0,
\qquad
\hbox{\rm $\P^*$-a.s.}
$$

Define the set of the minimizers of $U(\, \cdot\, )$:
\begin{equation}
    \mathscr{U}_{\min}
    :=
    \Big\{ x\in \T: \, U(x) = \min_{y\in \T} U(y) \Big\} \, .
    \label{Umin}
\end{equation}

\noindent Since $\inf_{x: \, |x|=n} U(x) \to \infty$ $\P^*$-a.s.\ (see (\ref{U->infty})), the set $\mathscr{U}_{\min}$ is finite and non-empty.

The main result of the paper is as follows.

\medskip

\begin{theorem}
\label{t:main}

 Assume $(\ref{cond-hab})$ and $(\ref{integrability-assumption})$.
 For any $\varepsilon>0$,\footnote{By convergence in $\P^*$-probability, we mean convergence in probability under $\P^*$.}
 $$
 \sup_{x\in \T}
 P_\omega 
 \Big\{ \, \Big| \frac{L_n(x)}{\frac{n}{\log n}} - \frac{\sigma^2}{4 D_\infty}\, \ee^{-U(x)} \Big| > \varepsilon \Big\}
 \to 
 0, 
 \qquad\hbox{in $\P^*$-probability} \, ,
 $$
 
 \noindent where $U(\, \cdot \, )$ is the symmetrized potential in $(\ref{U})$, $D_\infty$ the $\P^*$-almost sure positive limit of the derivative martingale $(D_n)$ in $(\ref{Dn})$, and
 \begin{equation}
     \sigma^2
     :=
     \E \Big( \sum_{y: \, |y|=1} V(y)^2 \, \ee^{-V(y)} \Big) \in (0, \, \infty) \, .
     \label{sigma}
 \end{equation}

\end{theorem}

\medskip

\begin{corollary}
\label{c:main}

 Assume $(\ref{cond-hab})$ and $(\ref{integrability-assumption})$. 
 If $\mathscr{A}_n$ is the set of the most visited sites at time $n$ as in $(\ref{ensemble_points_favoris})$, then
 $$
     \p^* ( \mathscr{A}_n \subset \mathscr{U}_{\min} ) 
     \to 
     1 \, ,
 $$

 \noindent where $\mathscr{U}_{\min}$ is the set of the minimizers of $U(\, \cdot\, )$ in $(\ref{Umin})$.

\end{corollary}

\medskip

Our results are not as strong as they might look like. For example, Theorem \ref{t:main} does not claim that $P_\omega \{ \sup_{x\in \T} | \frac{L_n(x)}{\frac{n}{\log n}} - \frac{\sigma^2}{4 D_\infty}\, \ee^{-U(x)} | > \varepsilon \} \to 0$ in $\P^*$-probability. It essentially says, in view of Proposition \ref{p:negligeable} below, that for any {\bf fixed} $x\in \T$, $P_\omega \{ | \frac{L_n(x)}{\frac{n}{\log n}} - \frac{\sigma^2}{4 D_\infty}\, \ee^{-U(x)} | > \varepsilon\} \to 0$ in $\P^*$-probability. Corollary \ref{c:main} is much weaker than what T\'oth~\cite{toth} proved for the symmetric Bernoulli random walk on $\z$: for example, it does not claim that $\p^*$-a.s., $\mathscr{A}_n \subset \mathscr{U}_{\min}$ for all sufficiently large $n$; we even do not know whether this is true.

For local time at fixed site of biased random walks on Galton--Watson trees in other recurrent regimes, see the recent paper \cite{ylocaltime}.

An important ingredient in the proof of Theorem \ref{t:main} is the following estimate on the local time of vertices that are away from the root:

\medskip

\begin{proposition}
\label{p:negligeable}

 Assume $(\ref{cond-hab})$ and $(\ref{integrability-assumption})$. 
 Then
 $$
    \lim_{\varepsilon\to 0}
    \limsup_{n\to \infty} \, 
    \p^* \Big\{ 
    \max_{x\in \T: \, U(x) \ge \log (\frac{8}{\varepsilon^2})}
    L_n (x)
    \ge \frac{\varepsilon \, n}{\log n} \Big\}
    =
    0 \, .
 $$

\end{proposition}

\medskip
 
Proposition \ref{p:negligeable} is given in Section \ref{s:preliminaries:BRW}. Theorem \ref{t:main} and Corollary \ref{c:main} are proved in Section \ref{s:pf}.
   
Throughout the paper, for any pair of vertices $x$ and $y$, we write $x<y$ or $y>x$ if $y$ is a (strict) descendant of $x$, and $x\le y$ or $y\ge x$ if either $y$ is either a (strict) descendant of $x$, or $x$ itself. For any $x\in\T$, we use $x_i$ (for $0\le i\le |x|$) to denote the ancestor of $x$ in the $i$-th generation; in particular, 
$x_0 = \varnothing$ and $x_{|x|} =x$.

\section{Proof of Proposition \ref{p:negligeable}}
\label{s:preliminaries:BRW}


$\phantom{aob}$ We start with some preliminaries. Define
\begin{equation}
    \Lambda(x)
    :=
    \sum_{y: \, {\buildrel \leftarrow \over y} = x} \ee^{-[V(y)-V(x)]} \, ,
    \qquad
    x\in \T \, ,
    \label{Lambda}
\end{equation}

\noindent In particular, $\Lambda(\varnothing) = \sum_{x: \, |x|=1} \ee^{-V(x)}$.

Let $S_i-S_{i-1}$, $i\ge 1$, be i.i.d.\ random variables whose law is characterized by
\begin{equation}
    \E \Big[ h ( S_1) \Big]
    =
    \E \Big[ \sum_{x\in \T: \, |x|=1} \ee^{-V(x)} h(V(x)) \Big] \, ,
    \label{joint-law:(S,Lambda)}
\end{equation}

\noindent for any Borel function $h: \, \r \to \r_+$.   

The following fact, quoted from  \cite{yzlocaltree}, is a variant of  the so-called ``many-to-one formula" for the branching random walk. 

\begin{fact}
\label{f:many-to-one-appli1}

 Assume $(\ref{cond-hab})$ and $(\ref{integrability-assumption})$. 
 Let $\Lambda(x)$ be as in $(\ref{Lambda})$.
 For any $n\ge 1$ and any Borel function $g: \r^{n+1} \to \r_+$, 
 we have
 $$
 \E \Big[ \sum_{x\in \T: \, |x|=n} g\Big( V(x_1), \, \cdots, \, V(x_n), \, \Lambda(x) \Big) \Big]
 =
 \E \Big[ \ee^{S_n} \, G\Big( S_1, \, \cdots, \, S_n\Big) \Big] \, ,
 $$
 where $S_i-S_{i-1}$, $i\ge 1$, are i.i.d.\ whose common distribution is given in $(\ref{joint-law:(S,Lambda)})$, and 
 $$
 G(a_1, \, \cdots, \, a_n) 
 := 
 \E  [ g(a_1, \, \cdots, \, a_n, \, \sum_{x\in \T: \, |x|=1} \ee^{-V(x)}) ] \, .
 $$

\end{fact}

\medskip

 Define  a reflecting barrier at (notation: $]\! ]\varnothing, \, x[\![ \; := \; ]\! ]\varnothing, \, x]\!] \backslash \{ x\}$)
\begin{eqnarray}
    \mathscr{L}_n^{(\gamma)}
 &:=& \Big\{ x: \, 
    \sum_{z\in \, ]\! ]  \varnothing, \, x]\! ]} \ee^{V(z)-V(x)} > \frac{n}{(\log n)^\gamma} , \;
    \nonumber
    \\
 &&\qquad\qquad 
    \sum_{z\in \, ]\! ] \varnothing, \, y]\! ]} \ee^{V(z)-V(y)} \le \frac{n}{(\log n)^\gamma}, \; 
    \forall y\in \, ]\! ]\varnothing, \, x[\![ \Big\} ,
    \label{gamma}
\end{eqnarray}

\noindent where $\gamma\in \r$ is a fixed parameter. We write $x< \mathscr{L}_n^{(\gamma)}$ if $\sum_{z\in \, ]\! ] \varnothing, \, y]\! ]} \ee^{V(z)-V(y)} \le \frac{n}{(\log n)^\gamma}$ for all $y\in \, ]\! ]\varnothing, \, x[\![\,$.

We recall two results from \cite{yzlocaltree}. The first justifies the presence of the barrier $\mathscr{L}_n^{(\gamma)}$ for the biased walk $(X_n)$, and the second describes the local time at the root.

\medskip

\begin{fact}[\cite{yzlocaltree}]
 \label{f:ligne-darret}
 
 Assume $(\ref{cond-hab})$ and $(\ref{integrability-assumption})$. 
If $\gamma<2$, then   
 $$
     \lim_{n\to \infty} \, 
     \p \Big( \bigcup_{i=1}^n\{ X_i \in \mathscr{L}^{(\gamma)}_n \} \Big)
     =
     0 \, .
 $$

\end{fact}

\medskip

\begin{fact}[\cite{yzlocaltree}]
\label{f:x=0}

 Assume $(\ref{cond-hab})$ and $(\ref{integrability-assumption})$. 
 For any $\varepsilon>0$,
 $$
 P_\omega 
 \Big\{ \, \Big| \frac{L_n(\varnothing)}{\frac{n}{\log n}} - \frac{\sigma^2}{4 D_\infty} \, \ee^{-U(\varnothing)}\Big| > \varepsilon \Big\}
 \to 
 0, 
 \qquad\hbox{in $\P^*$-probability} \, .
 $$

\end{fact}

\medskip

\begin{lemma}
\label{l:sum_iid}

 Let $0<a<1$ and $0<p<1$. Let $(\xi_i, \, i\ge 1)$ be an i.i.d.\ sequence of random variables with $\P(\xi_1 =0) =1-a$ and $\P(\xi_1 \ge k) = a \, p^{k-1}$, $\forall k\ge 1$. 
 
 Let $0<\varepsilon<1$. If $1-p > \frac{8}{\varepsilon} \, a$, then
 $$
 \P\Big\{ \sum_{i=1}^n \xi_i \ge \lceil \varepsilon n \rceil\Big\}
    \le
    6na\, \ee^{- \frac{(1-p)\varepsilon n}{8}}\, .
 $$

\end{lemma}

\medskip

\noindent {\it Proof.} Let $s\in [1, \, \frac1p )$. Then $\E(s^{\xi_1}) = 1- a + \frac{a(1-p)s}{1-ps}$. So 
$$
\P\Big\{ \sum_{i=1}^n \xi_i \ge k\Big\}
\le
\frac{1}{s^k}\, \E\Big[ s^{\sum_{i=1}^n \xi_i} \, {\bf 1}_{\{ \sum_{i=1}^n \xi_i >0\} } \Big]
=
\frac{[\E(s^{\xi_1})]^n - [\P \{ \xi_1=0\}]^n}{s^k}\, .
$$

\noindent Observe that
\begin{eqnarray*}
    [\E(s^{\xi_1})]^n - [\P \{ \xi_1=0\}]^n
 &=&\Big( 1- a + \frac{a(1-p)s}{1-ps} \Big)^n - (1-a)^n
    \\
 &\le&n \frac{a(1-p)s}{1-ps} \, \Big( 1- a + \frac{a(1-p)s}{1-ps} \Big)^{n-1} \, ,
\end{eqnarray*}

\noindent where, in the last line, we used $x^n - y^n \le n(x-y)x^{n-1}$ (for $0\le y\le x$). Hence
\begin{equation}
    \P\Big\{ \sum_{i=1}^n \xi_i \ge k\Big\}
    \le
    s^{-k} \, n \frac{a(1-p)s}{1-ps} \, \Big( 1- a + \frac{a(1-p)s}{1-ps} \Big)^{n-1} \, .
    \label{sum_iid-proof1}
\end{equation}

{\it First case:} $\frac13 \le p<1$. We take $s:= \frac{1+p}{2p}\in [1, \, \frac1p )$, so that $\frac{(1-p)s}{1-ps} = \frac{1+p}{p}$; hence by (\ref{sum_iid-proof1}),
\begin{eqnarray*}
    \P\Big\{ \sum_{i=1}^n \xi_i \ge k\Big\}
 &\le& \Big(\frac{1+p}{2p}\Big)^{\!-k} \; n \frac{a(1+p)}{p} \, \Big( 1+ \frac{a}{p} \Big)^{n-1} 
    \\
 &=& na\frac{1+p}{p} \Big( 1+\frac{1-p}{2p}\Big)^{\!-k} \, \Big( 1+ \frac{a}{p} \Big)^{n-1} 
    \\
 &\le& \frac{2na}{p} \Big( 1+\frac{1-p}{2p}\Big)^{\!-k} \, \Big( 1+ \frac{a}{p} \Big)^n\, .
\end{eqnarray*}

\noindent Since $(1+u)^{-1} \le \ee^{-u/2}$ (for $0\le u\le 1$) and $1+v \le \ee^v$ (for $v\ge 0$), applied to $u := \frac{1-p}{2p} \le 1$ and $v:= \frac{a}{p}$, we obtain, in case $\frac13 \le p<1$,
$$
    \P\Big\{ \sum_{i=1}^n \xi_i \ge k\Big\}
    \le
    6na\, \exp \Big( - \frac{(1-p)k}{4p} + \frac{na}{p} \Big) \, .
$$

{\it Second and last case:} $0<p\le \frac13$. We choose $s:= 2<\frac1p$, so $\frac{(1-p)s}{1-ps} \le 4$; by (\ref{sum_iid-proof1}), we obtain:
$$
\P\Big\{ \sum_{i=1}^n \xi_i \ge k\Big\}
\le
4na \, 2^{-k} (1+3a)^{n-1}
\le
4na \, 2^{-k} (1+3a)^n \, .
$$

\noindent In view of the inequality $1+v \le \ee^v$ (for $v\ge 0$; applied to $v:= 3a$), we obtain, in case $0<p\le \frac13$,
$$
    \P\Big\{ \sum_{i=1}^n \xi_i \ge k\Big\}
    \le
    4na \, 2^{-k} \, \ee^{3an}\, .
$$

So in both situations, as long as $1-p > \frac{8}{\varepsilon} \, a$, we have, for $k := \lceil \varepsilon n \rceil$, $- \frac{(1-p)k}{4p} + \frac{na}{p} \le -n( \frac{(1-p)\varepsilon}{4p} - \frac{a}{p} ) \le - \frac{(1-p)\varepsilon \, n}{8p} \le - \frac{(1-p)\varepsilon n}{8}$, and $2^{-k} \ee^{3an} \le \ee^{-n(\varepsilon \log 2 - 3a)} \le \ee^{-n(\varepsilon \log 2 - \frac{3\varepsilon}{8})}$, which is bounded by $\ee^{- \frac{\varepsilon n}{8}}$ (because $\log 2 \ge \frac12$), and a fortiori by $\ee^{- \frac{(1-p)\varepsilon n}{8}}$. Lemma \ref{l:sum_iid} is proved.\hfill$\Box$

\bigskip

We now proceed to the proof of Proposition \ref{p:negligeable}. Define
\begin{eqnarray}
    T_x
 &:=& \inf\{ i\ge 0: \, X_i =x \} \, , 
    \qquad x\in \T \, ,
    \label{T}
    \\
    T_\varnothing^+
 &:=& \inf\{ i\ge 1: \, X_i = \varnothing \} \, .
    \label{tau}
\end{eqnarray}

\noindent In words, $T_x$ is the first hitting time at $x$ by the biased walk, whereas $T_\varnothing^+$ is the first {\it return} time to the root $\varnothing$. 

Let $x\in \T \backslash\{ \varnothing\}$. The probability $P_\omega (T_x < T_\varnothing^+)$ only involves a one-dimensional random walk in random environment (namely, the restriction at $[\! [ \varnothing, \, x[\![\,$ of the biased walk $(X_i)$), so a standard result for one-dimensional random walks in random environment (Golosov~\cite{golosov}) tells us that
\begin{eqnarray}
    P_\omega (T_x < T_\varnothing^+)
 &=&\frac{\omega(\varnothing, \, x_1)\, \ee^{V(x_1)}}{\sum_{z\in \, ]\!] \varnothing, \, x]\!]} \ee^{V(z)}}    =
    \frac{\omega(\varnothing, \, {\buildrel \leftarrow \over \varnothing})}{\sum_{z\in \, ]\!] \varnothing, \, x]\!]} \ee^{V(z)}} \, ,
    \label{1D-MAMA}
    \\
    P_{x,\omega} \{ T_\varnothing < T_x^+ \}
 &=& \frac{\ee^{U(x)}}{\sum_{z\in \, ] \! ] \varnothing, \, x]\! ]} \ee^{V(z)}} \, ,
    \label{1D-MAMA-bis}
\end{eqnarray}
    
\noindent where $x_1$ is the ancestor of $x$ in the first generation.

\medskip

\noindent {\it Proof of Proposition \ref{p:negligeable}.} By Fact \ref{f:ligne-darret}, for all $\gamma_1<2$, we have $\p^* ( \cup_{i=1}^m\{ X_i \in \mathscr{L}_m^{(\gamma_1)}\} ) \to 0$, $m\to \infty$. So it suffices to check that for some $\gamma_1<2$,
$$
    \lim_{b\to 0}
    \limsup_{m\to \infty} \, 
    \p^* \Big\{ 
    \max_{x<\mathscr{L}_m^{(\gamma_1)}: \, U(x) \ge \log (\frac{8}{b^2})} L_m (x)
\ge \frac{bm}{\log m} \Big\}
    =
    0 \, .
$$

\noindent Since $\P^*(U(\varnothing) \ge \log (\frac{8}{b^2})) \to 0$ for $b\to 0$, it suffices to prove that for some $\gamma_1<2$,
\begin{equation}
    \lim_{b\to 0}
    \limsup_{m\to \infty} \, 
    \p^* \Big\{ 
    \max_{x\in\T \backslash \{ \varnothing\} : \, x<\mathscr{L}_m^{(\gamma_1)}, \, U(x) \ge \log (\frac{8}{b^2})} L_m (x)
\ge \frac{bm}{\log m} \Big\}
    =
    0 \, .
\label{pf-prop-negligeable-1}
\end{equation}

Let $T_\varnothing^{(0)} := 0$ and inductively $T_\varnothing^{(j)} := \inf\{ i> T_\varnothing^{(j-1)}: \, X_i = \varnothing \}$, for $j\ge 1$. In words, $T_\varnothing^{(j)}$ is the $j$-th return time to $\varnothing$. We have, for $n\ge 2$, $c>0$, $\varepsilon\in (0, \, 1)$, $1<\gamma<2$ and $m(n) = \lfloor c\, n \log n \rfloor$,
\begin{eqnarray*}
 &&\p^* \Big\{ 
    \max_{x\in\T \backslash \{ \varnothing\} : \, x<\mathscr{L}_{m(n)}^{(\gamma)}, \, U(x) \ge \log (\frac{8}{\varepsilon})} L_{m(n)} (x)
\ge \varepsilon n \Big\}
    \\
 &\le& \p^* \{ T_\varnothing^{(n)} \le m(n) \}
    +
    \p^* \Big\{ 
    \max_{x\in\T \backslash \{ \varnothing\} : \, x<\mathscr{L}_{m(n)}^{(\gamma)}, \, U(x) \ge \log (\frac{8}{\varepsilon})} L_{T_\varnothing^{(n)}} (x) \ge \varepsilon n \Big\} \, .
\end{eqnarray*}

\noindent By Fact \ref{f:x=0}, $\frac{T_\varnothing^{(n)}}{n\log n} \to \frac{4D_\infty}{\sigma^2} \ee^{U(\varnothing)}$ in $\p^*$-probability, so the portmanteau theorem implies that $\limsup_{n\to \infty} \p^*\{ T_\varnothing^{(n)} \le m(n) \} \le \P^*\{ \frac{4D_\infty}{\sigma^2} \ee^{U(\varnothing)} \le c\}$. Assume, for the time being, that we are able to prove that for some $\gamma<2$, any $c>0$ and any $0<\varepsilon<1$,
\begin{equation}
    \p^* \Big\{ 
    \max_{x\in\T \backslash \{ \varnothing\} : \, x<\mathscr{L}_{m(n)}^{(\gamma)}, \; U(x) \ge \log (\frac{8}{\varepsilon})}
    L_{T_\varnothing^{(n)}} (x)
    \ge \varepsilon n \Big\}
    \to 0 ,
    \qquad
    n\to \infty \, .
    \label{pf-prop-negligeable-2}
\end{equation}

\noindent Then we will have
$$
\limsup_{n\to \infty} 
\p^* \Big\{ 
\max_{x\in\T \backslash \{ \varnothing\} : \, x<\mathscr{L}_{m(n)}^{(\gamma)}, \, U(x) \ge \log (\frac{8}{\varepsilon})} L_{m(n)} (x)
\ge \varepsilon n \Big\}
\le
\P^*\Big\{ \frac{4D_\infty}{\sigma^2} \ee^{U(\varnothing)} \le c\Big\} \, .
$$

\noindent Since $n \le \frac2c \frac{m(n+1)}{\log m(n)}$ (for all sufficiently large $n$), this will yield
$$
\limsup_{n\to \infty} 
\p^* \Big\{ 
\max_{x\in\T \backslash \{ \varnothing\} : \, x<\mathscr{L}_{m(n)}^{(\gamma)}, \, U(x) \ge \log (\frac{8}{\varepsilon})} L_{m(n)} (x)
\ge \frac{2\varepsilon}{c} \frac{m(n+1)}{\log m(n)} \Big\}
\le
\P^*\Big\{ \frac{4D_\infty}{\sigma^2} \ee^{U(\varnothing)} \le c\Big\} \, .
$$

\noindent Let $m\in [m(n), \, m(n+1)] \cap \z$. Then $L_{m(n)} (x) \le L_m(x)$ (for all $x\in \T$); on the other hand, if $x<\mathscr{L}_{m}^{(\gamma_1)}$, then $x<\mathscr{L}_{m(n)}^{(\gamma)}$ for all $\gamma_1\in (\gamma, \, 2)$ and all sufficiently large $n$. Consequently, we will have, for all $c>0$ and $\varepsilon\in (0, \, 1)$,
$$
\limsup_{m\to \infty} 
\p^* \Big\{ 
\max_{x\in\T \backslash \{ \varnothing\} : \, x<\mathscr{L}_{m}^{(\gamma_1)}, \, U(x) \ge \log (\frac{8}{\varepsilon})} L_m (x)
\ge \frac{2\varepsilon}{c} \frac{m}{\log m} \Big\}
\le
\P^*\Big\{ \frac{4D_\infty}{\sigma^2} \ee^{U(\varnothing)} \le c\Big\} \, .
$$

\noindent Taking $c:= 2\varepsilon^{1/2}$ will then yield (\ref{pf-prop-negligeable-1}) (writing $b:= \varepsilon^{1/2}$ there) and thus Proposition \ref{p:negligeable}.

The rest of the section is devoted to the proof of (\ref{pf-prop-negligeable-2}). By \eqref{ascv}, $\frac{1}{(\log n)^3} \max_{0\le i\le n} |X_i|$ converges $\p^*$-a.s.\ to a positive constant, and since $\frac{T_\varnothing^{(n)}}{n\log n}$ converges in $\p^*$-probability to a positive limit, we deduce that $\frac{1}{(\log n)^3} \max_{0\le i\le T_\varnothing^{(n)}} |X_i|$ converges in $\p^*$-probability to a positive limit. So the proof of (\ref{pf-prop-negligeable-2}) is reduced to showing the following estimate: for some $1<\gamma<2$, any $c>0$ and any $0<\varepsilon<1$,
$$
    \p^* \Big\{ 
    \max_{x<\mathscr{L}_{m(n)}^{(\gamma)}: \; U(x) \ge \log (\frac{8}{\varepsilon}), \; 1\le |x| \le (\log n)^4}
    L_{T_\varnothing^{(n)}} (x)
    \ge \varepsilon n \Big\}
    \to 0 ,
    \qquad
    n\to \infty \, .
$$

For $k\ge 1$, we have
\begin{eqnarray*}
 &&P_\omega \Big\{ 
    \max_{x<\mathscr{L}_{m(n)}^{(\gamma)}: \; U(x) \ge \log (\frac{8}{\varepsilon}), \; 1\le |x| \le (\log n)^4}
    L_{T_\varnothing^{(n)}} (x) \ge k \Big\}
    \\
 &\le& \sum_{x<\mathscr{L}_{m(n)}^{(\gamma)}: \; U(x) \ge \log (\frac{8}{\varepsilon}), \; 1\le |x| \le (\log n)^4}
    P_\omega \{ L_{T_\varnothing^{(n)}} (x) \ge k\} \, .
\end{eqnarray*}

\noindent The law of $L_{T_\varnothing^{(n)}} (x)$ under $P_\omega$ is the law of $\sum_{i=1}^n \xi_i$, where $(\xi_i, \, i\ge 1)$ is an i.i.d.\ sequence with $P_\omega(\xi_1 =0) =1-a$ and $P_\omega(\xi_1 \ge k) = a \, p^{k-1}$, $\forall k\ge 1$, where 
\begin{eqnarray*}
    1-p 
 &:=& P_{x,\omega} \{ T_\varnothing < T_x^+ \}
    =
    \frac{\ee^{U(x)}}{\sum_{z\in \, ] \! ] \varnothing, \, x]\! ]} \ee^{V(z)}} \, ,
    \\
    a
 &:=& P_\omega \{ T_x < T_\varnothing^+\} 
    =
    \frac{\omega(\varnothing, \, {\buildrel \leftarrow \over \varnothing})}{\sum_{z\in \, ] \! ] \varnothing, \, x]\! ]} \ee^{V(z)}} \, .
\end{eqnarray*}

\noindent [We have used \eqref{1D-MAMA} and \eqref{1D-MAMA-bis}.]
 
If $U(x) \ge \log (\frac{8}{\varepsilon})$, then $1-p > \frac{8}{\varepsilon} \, a$, so we are entitled to apply Lemma \ref{l:sum_iid} to arrive at:
\begin{eqnarray*}
 &&P_\omega \Big\{ 
    \max_{x<\mathscr{L}_{m(n)}^{(\gamma)}: \; U(x) \ge \log (\frac{8}{\varepsilon}), \; 1\le |x| \le (\log n)^4}
    L_{T_\varnothing^{(n)}} (x) \ge \lceil \varepsilon n\rceil \Big\}
    \\
 &\le&6n
    \sum_{x<\mathscr{L}_{m(n)}^{(\gamma)}: \; 1\le |x| \le (\log n)^4}
    \frac{\omega(\varnothing, \, {\buildrel \leftarrow \over \varnothing})}{\sum_{z\in \, ] \! ] \varnothing, \, x]\! ]} \ee^{V(z)}} 
    \exp\Big( - \frac{\varepsilon n}{8} \frac{\ee^{U(x)}}{\sum_{z\in \, ] \! ] \varnothing, \, x]\! ]} \ee^{V(z)}} \Big) \, .
\end{eqnarray*}

\noindent We have $\omega(\varnothing, \, {\buildrel \leftarrow \over \varnothing}) \le 1$. It remains to check the following convergence in $\P^*$-probability (for $n\to \infty$):
\begin{equation}
    \sum_{x<\mathscr{L}_{m(n)}^{(\gamma)}: \; 1\le |x| \le (\log n)^4}
    \frac{n}{\sum_{z\in \, ] \! ] \varnothing, \, x]\! ]} \ee^{V(z)}} 
    \exp\Big( - \frac{\varepsilon}{8} \frac{n\, \ee^{U(x)}}{\sum_{z\in \, ] \! ] \varnothing, \, x]\! ]} \ee^{V(z)}} \Big)
    \to
    0 \, .
    \label{pf-prop-negligeable-3}
\end{equation}

\noindent Recall the definition of $\mathscr{L}_{m(n)}^{(\gamma)}$: $x<\mathscr{L}_{m(n)}^{(\gamma)}$ implies $\frac{\ee^{V(x)}}{\sum_{z\in \, ]\! ] \varnothing, \, x]\! ]} \ee^{V(z)}} \ge \frac{(\log m(n))^\gamma}{m(n)}$, which is $\ge \frac{(\log n)^{\gamma-1}}{cn}$ for all sufficiently large $n$ (say $n\ge n_0$). Also, we recall that $\ee^{U(x)} = \frac{\ee^{V(x)}}{1+\Lambda(x)}$, with $\Lambda(x) := \sum_{y: \, {\buildrel \leftarrow \over y} =x} \ee^{-[V(y)-V(x)]}$ as in (\ref{Lambda}).

For the sum $\sum_{x<\mathscr{L}_{m(n)}^{(\gamma)}}$ on the left-hand side of (\ref{pf-prop-negligeable-3}), we distinguish two possible situations depending on the value of $\Lambda(x)$. Let $0<\varrho<1$. Applying the elementary inequality $\lambda \ee^{-\lambda} \le c_2\, \ee^{-\lambda/2}$ (for $\lambda\ge 0$) to $\lambda := \frac{\varepsilon}{8} \, \frac{n\, \ee^{U(x)}}{\sum_{z\in \, ] \! ] \varnothing, \, x]\! ]} \ee^{V(z)}}$, we see that for $n\ge n_0$,
\begin{eqnarray*}
 &&\sum_{x<\mathscr{L}_{m(n)}^{(\gamma)}} 
    {\bf 1}_{\{ 1+\Lambda(x) \le (\frac{n\, \ee^{V(x)}}{\sum_{z\in \, ] \! ] \varnothing, \, x]\! ]} \ee^{V(z)}})^{\varrho}\} } 
    \frac{n}{\sum_{z\in \, ] \! ] \varnothing, \, x]\! ]} \ee^{V(z)}} 
    \exp\Big( - \frac{\varepsilon}{8} \frac{n\, \ee^{U(x)}}{\sum_{z\in \, ] \! ] \varnothing, \, x]\! ]} \ee^{V(z)}} \Big)
    \\
 &\le& c_2 \, \sum_{x<\mathscr{L}_{m(n)}^{(\gamma)}} 
    {\bf 1}_{\{ 1+\Lambda(x) \le (\frac{n\, \ee^{V(x)}}{\sum_{z\in \, ] \! ] \varnothing, \, x]\! ]} \ee^{V(z)}})^{\varrho}\} } \,
    \frac{8}{\varepsilon} \, \ee^{-U(x)} 
    \exp\Big( - \frac{\varepsilon}{16(1+\Lambda(x))} \frac{n\, \ee^{V(x)}}{\sum_{z\in \, ] \! ] \varnothing, \, x]\! ]} \ee^{V(z)}} \Big)
    \\
 &\le& \frac{8c_2}{\varepsilon} \, \sum_{x<\mathscr{L}_{m(n)}^{(\gamma)}} 
    \ee^{-U(x)} 
    \exp\Big( - \frac{\varepsilon}{16} (\frac{n\, \ee^{V(x)}}{\sum_{z\in \, ] \! ] \varnothing, \, x]\! ]} \ee^{V(z)}})^{1-\varrho} \Big) \, .
\end{eqnarray*}

\noindent Since $\frac{n\, \ee^{V(x)}}{\sum_{z\in \, ] \! ] \varnothing, \, x]\! ]} \ee^{V(z)}} \ge \frac1c (\log n)^{\gamma-1}$ (for $x<\mathscr{L}_{m(n)}^{(\gamma)}$ and $n\ge n_0$), this yields, for $n\ge n_0$,
\begin{eqnarray*}
 &&\sum_{x<\mathscr{L}_{m(n)}^{(\gamma)}} 
    {\bf 1}_{\{ 1+\Lambda(x) \le (\frac{n\, \ee^{V(x)}}{\sum_{z\in \, ] \! ] \varnothing, \, x]\! ]} \ee^{V(z)}})^{\varrho}\} } 
    \frac{n}{\sum_{z\in \, ] \! ] \varnothing, \, x]\! ]} \ee^{V(z)}} 
    \exp\Big( - \frac{\varepsilon}{8} \frac{n\, \ee^{U(x)}}{\sum_{z\in \, ] \! ] \varnothing, \, x]\! ]} \ee^{V(z)}} \Big)
    \\
 &\le& \frac{8c_2}{\varepsilon} \, \exp\Big( - \frac{\varepsilon}{16c^{1-\varrho}} (\log n)^{(\gamma-1)(1-\varrho)} \Big)
    \sum_{x<\mathscr{L}_{m(n)}^{(\gamma)}} \ee^{-U(x)} \, ,
\end{eqnarray*}

\noindent which converges to $0$ in $\P^*$-probability (recalling that for any $\gamma\in \r$, $\frac{1}{\log n} \sum_{x\in \T: \, x < \mathscr{L}_n^{(\gamma)}} \ee^{-U(x)}$ converges in $\P^*$-probability to a finite limit; see \cite{yzlocaltree}). So it remains to prove that there exists $\varrho \in (0, \, 1)$ such that (removing the big exponential term which is bounded by $1$)
$$
\sum_{x<\mathscr{L}_{m(n)}^{(\gamma)}: \; 1\le |x| \le (\log n)^4} 
{\bf 1}_{\{ 1+\Lambda(x) > (\frac{n\, \ee^{V(x)}}{\sum_{z\in \, ] \! ] \varnothing, \, x]\! ]} \ee^{V(z)}})^{\varrho}\} } 
\frac{n}{\sum_{z\in \, ] \! ] \varnothing, \, x]\! ]} \ee^{V(z)}} 
\to
0 \, ,
$$

\noindent in $\P^*$-probability (for $n\to \infty$). Since $\lim_{r\to \infty} \inf_{|x|=r} V(x) \to \infty$ $\P^*$-a.s.\ (see \eqref{U->infty}), it suffices to prove the existence of $\varrho \in (0, \, 1)$ and $\gamma\in (1, \, 2)$ such that for all $\alpha>0$ and $n\to \infty$,
\begin{equation}
    \sum_{x<\mathscr{L}_{m(n)}^{(\gamma)}: \; 1\le |x| \le (\log n)^4}
    {\bf 1}_{\{ 1+\Lambda(x) > (\frac{n\, \ee^{V(x)}}{\sum_{z\in \, ] \! ] \varnothing, \, x]\! ]} \ee^{V(z)}})^{\varrho}\} } 
    \frac{n}{\sum_{z\in \, ] \! ] \varnothing, \, x]\! ]} 
    \ee^{V(z)}} \,
    {\bf 1}_{\{ \underline{V}(x) \ge -\alpha\} } 
    \to
    0 \, ,
    \label{proof1}
\end{equation}

\noindent in $\P^*$-probability, where $\underline{V}(x) := \min_{z\in \, ] \! ] \varnothing, \, x]\!]} V(z)$.

To prove this, we first recall that $x<\mathscr{L}_{m(n)}^{(\gamma)}$ implies that for all $y\in \, ]\!] \varnothing, \, x]\!]$, we have $\frac{\sum_{z\in \, ]\! ] \varnothing, \, y]\! ]} \ee^{V(z)}}{\ee^{V(y)}} \le \frac{cn}{(\log n)^{\gamma-1}}$ (for $n\ge n_0$) which is bounded by $n$ for all sufficiently large $n$ (say $n\ge n_1$); a fortiori $\overline{V}(y) - V(y) \le \log n$ (with $\overline{V}(y):= \max_{z\in \, ]\!]  \varnothing, \, y]\! ]} V(y)$). By Fact \ref{f:many-to-one-appli1}, we obtain, for $n\ge n_0\vee n_1$,
\begin{eqnarray*}
 &&\E \Big[ \sum_{x<\mathscr{L}_{m(n)}^{(\gamma)}, \; 1\le |x| \le (\log n)^4}
    {\bf 1}_{\{ 1+\Lambda(x) > (\frac{n\, \ee^{V(x)}}{\sum_{z\in \, ] \! ] \varnothing, \, x]\! ]} \ee^{V(z)}})^{\varrho}\} } 
    \frac{n}{\sum_{z\in \, ] \! ] \varnothing, \, x]\! ]} \ee^{V(z)}} \,
    {\bf 1}_{\{ \underline{V}(x) \ge -\alpha\} }
    \Big]
    \\
 &\le& \sum_{k=1}^{\lfloor (\log n)^4\rfloor}
    \E \Big[ 
    \sum_{x: \, |x|=k} 
    {\bf 1}_{ \{ \overline{V}(y) - V(y) \le \log n, \; \forall y\in \, ]\!] \varnothing, \, x]\!] \} } \times
    \\
 && \qquad {\bf 1}_{\{ 1+\Lambda(x) > (\frac{n\, \ee^{V(x)}}{\sum_{z\in \, ] \! ] \varnothing, \, x]\! ]} \ee^{V(z)}})^{\varrho}\} } \,
    \frac{n}{\sum_{z\in \, ] \! ] \varnothing, \, x]\! ]} \ee^{V(z)}} \,
    {\bf 1}_{\{ \underline{V}(x) \ge -\alpha\} }
    \Big]
    \\
 &=&\sum_{k=1}^{\lfloor (\log n)^4\rfloor}
    \E \Big[ 
    \ee^{S_k} \,
    {\bf 1}_{ \{ S^{\#}_k \le \log n\} } \,
    F\Big( (\frac{n\, \ee^{S_k}}{\sum_{i=1}^k \ee^{S_i}})^{\varrho} \Big)
\frac{n}{\sum_{i=1}^k \ee^{S_i}} \,
    {\bf 1}_{\{ \underline{S}_k \ge -\alpha\} }
    \Big] \, ,
\end{eqnarray*}

\noindent where $F(\lambda) := \P( 1+ \sum_{x: \, |x|=1}\ee^{-V(x)} > \lambda)$ for $\lambda >0$,  $\overline{S}_k := \max_{1\le i\le k}S_i$, $\underline{S}_k := \max_{1\le i\le k}S_i$, and $S^{\#}_k:=\max_{1\le i \le k} (\overline S_i- S_i)$ for any $k\ge 1$. 

An application of the H\"older inequality, using assumption (\ref{integrability-assumption}), yields the existence of $\delta_1>0$ such that
\begin{equation}
    \E\Big[ \Big( \sum_{x: \, |x|=1} \ee^{-V(x)} \Big)^{1+\delta_1} \Big]
    <
    \infty\, .
    \label{integrabilite_sous_Q}
\end{equation}

\noindent As such, $c_3 := \E [ ( 1+\sum_{x: \, |x|=1} \ee^{-V(x)} )^{1+\delta_1} ] <\infty$, so $F(\lambda) \le c_3 \, \lambda^{-1-\delta_1}$ for all $\lambda>0$. Consequently,
\begin{eqnarray}
 &&\E \Big[ \sum_{x<\mathscr{L}_{m(n)}^{(\gamma)}, \; 1\le |x| \le (\log n)^4}
    {\bf 1}_{\{ 1+\Lambda(x) > (\frac{n\, \ee^{V(x)}}{\sum_{z\in \, ] \! ] \varnothing, \, x]\! ]} \ee^{V(z)}})^{\varrho}\} } 
    \frac{n}{\sum_{z\in \, ] \! ] \varnothing, \, x]\! ]} \ee^{V(z)}} \,
    {\bf 1}_{\{ \underline{V}(x) \ge -\alpha\} }
    \Big]
    \nonumber
    \\
 &&\qquad \le c_3 \, \sum_{k=1}^{\lfloor (\log n)^4\rfloor}
    \E \Big[ 
    \Big( \frac{\sum_{i=1}^k \ee^{S_i}}{n\, \ee^{S_k}} \Big)^{\varrho(1+\delta_1)-1} \, 
    {\bf 1}_{ \{ S^{\#}_k \le \log n\} } \,
    {\bf 1}_{\{ \underline{S}_k \ge -\alpha\} } \Big] \, .
    \label{proof2}
\end{eqnarray}

\medskip

\begin{lemma}
\label{l:lemme_technique}

 Let $\delta$ be the constant in assumption $(\ref{integrability-assumption})$. For all $\alpha>0$ and $\delta_2 \in (0, \, \delta \wedge \frac{1}{16})$,
$$
    \lim_{n\to \infty} \,
    \sum_{k=1}^{\lfloor (\log n)^4\rfloor}
    \E \Big[ 
    \Big( \frac{\sum_{i=1}^k \ee^{S_i}}{n\, \ee^{S_k}} \Big)^{\delta_2} \, 
    {\bf 1}_{ \{ S^{\#}_k \le \log n\} } \,
    {\bf 1}_{\{ \underline{S}_k \ge -\alpha\} } \Big] 
    =
    0 \, .
$$

\end{lemma}

\medskip

Since it is possible to choose $0<\varrho<1$ such that $\varrho(1+\delta_1)-1$ lies in $(0, \, \delta \wedge \frac{1}{16})$, we can apply Lemma \ref{l:lemme_technique} to see that \eqref{proof2} implies \eqref{proof1}, and thus yields Proposition \ref{p:negligeable}.

It remains to prove Lemma \ref{l:lemme_technique}.

\bigskip

\noindent {\it Proof of Lemma \ref{l:lemme_technique}.} Since $\sum_{i=1}^k \ee^{S_i} \le k \, \overline{S}_k$, it suffices to check that
$$
\frac{(\log n)^{4\delta_2}}{n^{\delta_2}}
\sum_{k=1}^{\lfloor (\log n)^4\rfloor}
    \E \Big[ 
    \ee^{\delta_2(\overline{S}_k -S_k)} \, 
    {\bf 1}_{ \{ S^{\#}_k \le \log n\} } \,
    {\bf 1}_{\{ \underline{S}_k \ge -\alpha\} } \Big] 
\to
0 \, .
$$

Recall the law of $S_1$ from \eqref{joint-law:(S,Lambda)}. By assumption \eqref{integrability-assumption} and H\"older's inequality, we have 
$$
\E(\ee^{ a S_1})<\infty,
\qquad \forall a\in (-\delta, \, 1+\delta),
$$

\noindent where $\delta>0$ is the constant in \eqref{integrability-assumption}. In particular, $\E(\ee^{ a |S_1|})<\infty$ for all $0\le a< \delta$. Since $0<\delta_2<\delta$, we have $\E(\ee^{\delta_2(\overline{S}_k -S_k)}) \le \ee^{c_4 \, k}$ for some constant $c_4>0$ and all $k\ge 1$. So $\frac{(\log n)^{4\delta_2}}{n^{\delta_2}} \sum_{k=1}^{\lfloor (\log n)^{1/2}\rfloor} \E [ \ee^{\delta_2(\overline{S}_k -S_k)} ] \to 0$. It remains to prove that
$$
\frac{(\log n)^{4\delta_2}}{n^{\delta_2}}
\sum_{k=\lfloor (\log n)^{1/2}\rfloor}^{\lfloor (\log n)^4\rfloor}
    \E \Big[ 
    \ee^{\delta_2(\overline{S}_k -S_k)} \, 
    {\bf 1}_{ \{S^{\#}_k \le \log n\} } \,
    {\bf 1}_{\{ \underline{S}_k  \ge -\alpha\} } \Big] 
\to
0 \, .
$$

\noindent We make a change of indices $k=\lfloor (\log n)^{1/2}\rfloor + \ell$. Let $\widetilde{S}_\ell := S_{\ell + \lfloor (\log n)^{1/2}\rfloor} - S_{\lfloor (\log n)^{1/2}\rfloor}$, $\ell \ge 0$. Then $(\widetilde{S}_\ell, \, \ell \ge 0)$ is a random walk having the law of $(S_\ell, \, \ell \ge 0)$, and is independent of $(S_i, \, 1\le i\le \lfloor (\log n)^{1/2}\rfloor)$.     For $\ell \ge 0$, $\overline{S}_{\ell + \lfloor (\log n)^{1/2}\rfloor} -S_{\ell + \lfloor (\log n)^{1/2}\rfloor} =  \max(x,\,  \max_{0\le j\le \ell} \widetilde{S}_j) - \widetilde{S}_\ell \ge \max_{0\le j\le \ell} \widetilde{S}_j - \widetilde{S}_\ell$, where  $x:= \overline{S}_{\lfloor (\log n)^{1/2}\rfloor} -S_{\lfloor (\log n)^{1/2}\rfloor}$. So for $k\ge \lfloor (\log n)^{1/2}\rfloor$ and $\ell := k- \lfloor (\log n)^{1/2}\rfloor$, on the event that $\{ \underline{S}_k  \ge -\alpha\} $, either $ \max_{0\le j\le \ell} \widetilde{S}_j\le x$, then $\overline{S}_k -S_k=x-\widetilde{S}_\ell=\overline{S}_{\lfloor (\log n)^{1/2}\rfloor} -S_k\le \overline{S}_{\lfloor (\log n)^{1/2}\rfloor} +\alpha$, or $ \max_{0\le j\le \ell} \widetilde{S}_j> x$, then $\overline{S}_k -S_k=\max_{0\le j\le \ell} \widetilde{S}_j - \widetilde{S}_\ell$.  It follows that 
\begin{eqnarray*}
 &&\E \Big[ 
    \ee^{\delta_2(\overline{S}_k -S_k)} \, 
    {\bf 1}_{ \{ S^{\#}_k \le \log n\} } \,
    {\bf 1}_{\{ \underline{S}_{\lfloor (\log n)^{1/2}\rfloor} \ge -\alpha\} } \Big]
    \\
 &\le& \E (\ee^{\delta_2(\alpha+\overline{S}_{\lfloor (\log n)^{1/2}\rfloor})})
    +
    \P( \underline{S}_{\lfloor (\log n)^{1/2}\rfloor} \ge -\alpha) \times \E \Big[ 
    \ee^{\delta_2(\overline{S}_\ell -S_\ell)} \, 
    {\bf 1}_{ \{ S^{\#}_\ell \le \log n\} } \Big] \, .
\end{eqnarray*}

\noindent Since $\E (\ee^{\delta_2\overline{S}_{\lfloor (\log n)^{1/2}\rfloor}}) \le \ee^{c_4 \, (\log n)^{1/2}}$, we have $\frac{(\log n)^{4\delta_2}}{n^{\delta_2}} \sum_{k=\lfloor (\log n)^{1/2}\rfloor}^{\lfloor (\log n)^4\rfloor} \E (\ee^{\delta_2(\alpha+\overline{S}_{\lfloor (\log n)^{1/2}\rfloor})}) \to 0$. On the other hand, $\P( \underline{S}_{\lfloor (\log n)^{1/2}\rfloor} \ge -\alpha) \le c_5 \, (\log n)^{-1/4}$ for some constant $c_5>0$ and all $n\ge 2$ (see Kozlov~\cite{kozlov}); it suffices to prove that
$$
\frac{(\log n)^{4\delta_2-(1/4)}}{n^{\delta_2}}
\sum_{\ell=0}^\infty
    \E \Big[ 
    \ee^{\delta_2(\overline{S}_\ell -S_\ell)} \, 
    {\bf 1}_{ \{ S^{\#}_\ell \le \log n\} } \Big] 
\to
0 \, .
$$

\noindent This will be a straightforward consequence of the following estimate (applied to $\lambda := \log n$ and $b:= \delta_2$; it is here we use the condition $\delta_2 < \frac{1}{16}$): for any $0<b<\delta$,
\begin{equation}
    \limsup_{\lambda\to \infty} 
    \E \Big( \sum_{\ell =0}^{ \tau_\lambda-1}
    \ee^{-b[\lambda-(\overline{S}_\ell -S_\ell)]}
    \Big)
    <
    \infty \, ,
    \label{estimee_chute}
\end{equation}

\noindent where $\tau_\lambda := \inf\{ i\ge 1: \, \overline{S}_i - S_i > \lambda\}$. 

To prove (\ref{estimee_chute}), we define the (strictly) ascending ladder times $(H_i, i\ge0)$: $H_0:=0$ and for any $i\ge1$, $$ H_i:= \inf\{\ell>H_{i-1}: S_\ell> \max_{0\le j\le H_{i-1}} S_j\}.$$

\noindent Therefore,
$$
\E \Big( \sum_{\ell =0}^{ \tau_\lambda-1}
\ee^{ b (\overline{S}_\ell -S_\ell) }
\Big)
=
\sum_{i=1}^\infty  \E \Big( \sum_{\ell =H_{i-1}}^{H_i-1}    \ee^{ b (S_{H_{i-1}} -S_\ell) } 1_{\{S^{\#}_\ell\le \lambda\}}   \Big).
$$

\noindent We apply the strong Markov property, first at time $H_{i-1}$ to see that 
$$
\E \Big( \sum_{\ell =H_{i-1}}^{H_i-1} \ee^{ b (S_{H_{i-1}} -S_\ell) } 1_{\{S^{\#}_\ell\le \lambda\}}   \Big) 
\le 
\P\Big(S^{\#}_{H_{i-1}}\le \lambda\Big) \, \E \Big( \sum_{\ell =0}^{H_1-1}    \ee^{-b  S_\ell  } 1_{\{S^{\#}_\ell\le \lambda\}}   \Big) \, ,
$$

\noindent and then successively at times $H_1$, $H_2$, $\cdots$, $H_{i-1}$ to see that $\P (S^{\#}_{H_{i-1}}\le \lambda ) \le \P (S^{\#}_{H_1}\le \lambda )^{i-1}$. As such,
$$
\E \Big( \sum_{\ell =0}^{ \tau_\lambda-1}
\ee^{ b (\overline{S}_\ell -S_\ell) }
\Big)
\le
\sum_{i=1}^\infty 
\P\Big(S^{\#}_{H_1}\le \lambda\Big)^{i-1} \, 
\E \Big( \sum_{\ell =0}^{H_1-1}    \ee^{-b  S_\ell  } 1_{\{S^{\#}_\ell\le \lambda\}}   \Big) \, .
$$

\noindent We define $\sigma_{-\lambda}:= \inf\{n\ge0: S_n < - \lambda\}$. Then $1-\P (S^{\#}_{H_1}\le \lambda )=\P (\sigma_{-\lambda} < H_1 ) \ge {c_6\over 1+\lambda}$ for some constant $c_6>0$ and all sufficiently large $\lambda$, say $\lambda\ge \lambda_0$ (for the last elementary inequality, see for example, Lemma A.1 in \cite{yzenergy}). Thus we get that $$   \E \Big( \sum_{\ell =0}^{ \tau_\lambda-1} \ee^{ b (\overline{S}_\ell -S_\ell) } \Big) \le \frac{1+\lambda}{c_6}\;  \E \Big( \sum_{\ell =0}^{H_1-1}    \ee^{-b  S_\ell  } 1_{\{S^{\#}_\ell\le \lambda\}}   \Big).$$

\noindent Finally,  for all small $b>0$, there exists  some positive constant $c_7=c_7(b)>0$ such that
$$
\E \Big( \sum_{\ell =0}^{H_1-1} \ee^{-b S_\ell} {\bf 1}_{\{S^{\#}_\ell\le \lambda\}}   \Big) 
=
\E \Big( \sum_{\ell =0}^{H_1-1} \ee^{-b  S_\ell} {\bf 1}_{\{ \sigma_{-\lambda} > \ell\}} \Big) \le \frac{c_7}{\lambda}\, \ee^{b \lambda},
$$

\noindent by applying \cite{AHZ} (Lemma 6,  formula (4.17)) to $(-S_i, \, i\ge 1)$. This yields (\ref{estimee_chute}), and completes the proof of Lemma \ref{l:lemme_technique} and Proposition \ref{p:negligeable}.\hfill$\Box$

\section{Proof of Theorem \ref{t:main} and Corollary \ref{c:main}}
   \label{s:pf}


\noindent {\it Proof of Theorem \ref{t:main}.} Recall that $\lim_{k\to \infty} \inf_{x: \, |x|=k} U(x) \to \infty$ $\P^*$-a.s.\ (see \eqref{U->infty}). In view of Proposition \ref{p:negligeable}, we only need to prove that for any fixed $x\in \T$ and $\varepsilon>0$, when $n\to \infty$,
$$
P_\omega 
\Big\{ \, \Big| \frac{L_n(x)}{\frac{n}{\log n}} - \frac{\sigma^2}{4 D_\infty}\, \ee^{-U(x)} \Big| > \varepsilon \Big\}
\to 
0, 
\qquad\hbox{in $\P^*$-probability} \, .
$$

\noindent According to Fact \ref{f:x=0}, this is equivalent to convergence in $\P^*$-probability $P_\omega \{ \, | \frac{L_n(x)}{L_n(\varnothing)} - \ee^{-[U(x)-U(\varnothing)]} | > \varepsilon \} \to 0$ (for $n\to \infty$), and thus to the following statement: for any $x\in \T$ and $m\to \infty$,
$$
    P_\omega 
    \Big\{ \, \Big| \frac{L_{T_\varnothing^{(m)}}(x)}{m} - \ee^{-[U(x)-U(\varnothing)]} \Big| > \varepsilon \Big\}
    \to 
    0, 
    \qquad\hbox{in $\P^*$-probability} \, ,
$$

\noindent where $T_\varnothing^{(m)}$ is as before, the $m$-th return time of the biased walk $(X_i)$ to the root $\varnothing$. This, however, holds trivially as $L_{T_\varnothing^{(m)}}(x) - L_{T_\varnothing^{(m-1)}}(x)$, $m\ge 1$, are i.i.d.\ random variables under $P_\omega$ with $E_\omega[L_{T_\varnothing^{(1)}}(x)] = \ee^{-[U(x)-U(\varnothing)]}$. Theorem \ref{t:main} is proved.\hfill$\Box$

\bigskip

\noindent {\it Proof of Corollary \ref{c:main}.} Let $\varepsilon>0$ and $0<a<\frac12$. Let
$$
E_n(\varepsilon, \, a)
:=
\Big\{ \omega: \, \sup_{x\in \T} P_\omega \Big( \Big| \frac{L_n(x)}{\frac{n}{\log n}} - \frac{\sigma^2}{4 D_\infty}\, \ee^{-U(x)} \Big| > \varepsilon \Big) < a \Big\} \, .
$$

\noindent By Theorem \ref{t:main}, $\P^*(E_n(\varepsilon, \, a))\to 1$, $n\to \infty$. 

Let $x_n \in \mathscr{A}_n$, and let $x_{\min} \in \mathscr{U}_{\min}$. For all $\omega\in E_n(\varepsilon, \, a)$, we have
$$
P_\omega \Big( \Big| \frac{L_n(y)}{\frac{n}{\log n}} - \frac{\sigma^2}{4 D_\infty}\, \ee^{-U(y)} \Big| \le \varepsilon \Big)
\ge
1-a \, ,
$$

\noindent for $y= x_n$ and for $y= x_{\min}$; hence, for all $\omega\in E_n(\varepsilon, \, a)$,
$$
P_\omega \Big( \frac{L_n(x_n)}{\frac{n}{\log n}} \le \frac{\sigma^2}{4 D_\infty}\, \ee^{-U(x_n)} + \varepsilon, \; \frac{L_n(x_{\min})}{\frac{n}{\log n}} \ge \frac{\sigma^2}{4 D_\infty}\, \ee^{-U(x_{\min})} - \varepsilon \Big)
\ge
1-2a \, .
$$

\noindent By definition, $L_n(x_n) = \sup_{x\in \T} L_n(x) \ge L_n(x_{\min})$. Therefore, for all $\omega$,
$$
P_\omega \Big( \frac{\sigma^2}{4 D_\infty}\, \ee^{-U(x_n)} \ge \frac{\sigma^2}{4 D_\infty}\, \ee^{-U(x_{\min})} - 2\varepsilon \Big)
\ge
(1-2a)\, {\bf 1}_{E_n(\varepsilon, \, a)}(\omega) \, .
$$

\noindent Taking expectation with respect to $\P^*$ on both sides gives that
$$
\p^*\Big( \frac{\sigma^2}{4 D_\infty}\, \ee^{-U(x_n)} \ge \sup_{x\in \T} \frac{\sigma^2}{4 D_\infty}\, \ee^{-U(x)} - 2\varepsilon \Big)
\ge
(1-2a)\, \P^*(E_n(\varepsilon, \, a)) \, ,
$$

\noindent which converges to $1-2a$ when $n\to \infty$. Since $a>0$ can be as small as possible, this yields $\frac{\sigma^2}{4 D_\infty}\, \ee^{-U(x_n)} \to\sup_{x\in \T}  \frac{\sigma^2}{4 D_\infty}\, \ee^{-U(x)}$ in probability under $\p^*$, i.e., $U(x_n) \to \inf_{x\in \T} U(x)$ in probability under $\p^*$. 

Since $\mathscr{U}_{\min}$, the set of the minimizers of $U(\, \cdot \, )$, is $\P^*$-a.s.\ finite, we have $\inf_{x\in \mathscr{U}_{\min}} U(x) < \inf_{x\in \T \backslash \mathscr{U}_{\min}} U(x)$ $\P^*$-a.s., which yields $\p^*(x_n \in \mathscr{U}_{\min}) \to 1$, $n\to \infty$.\hfill$\Box$

\end{document}